\documentclass[12pt]{article}
\usepackage{amsmath,amssymb,amsthm,amscd}
\usepackage[small,nohug,heads=littlevee,noPostScript]{diagrams}
\title{Extremal contractions, stratified Mukai flops and
 Springer maps}
\author{Baohua Fu}
\newtheorem{Thm}{Theorem}[section]
\newtheorem{Lem}[Thm]{Lemma}
\newtheorem{Prop}[Thm]{Proposition}
\newtheorem{Cor}[Thm]{Corollary}

\newtheorem{Rque}[Thm]{Remark}

\newtheorem{Exam}[Thm]{Example}

\def\qit{{\mathbb Q}}
\def\zit{{\mathbb Z}}
\def\pit{{\mathbb P}}

\def\0{{\mathcal O}}
\def\g{{\mathfrak g}}
\def\h{{\mathfrak h}}
\def\p{{\mathfrak p}}
\def\l{{\mathfrak l}}
\def\n{{\mathfrak n}}

\def\L{{\mathcal L}}
\def\mtp{\mathop{\rm mtp}\nolimits}
\def\Ind{\mathop{\rm Ind}\nolimits}
\def\q{{\mathfrak q}}
\def\u{{\mathfrak u}}
\def\b{{\mathfrak b}}
\begin{document}
\maketitle
\begin{abstract}
We prove that two Springer maps of the same degree
 over a  nilpotent orbit closure are connected
by stratified Mukai flops, and the latter is 
 obtained by contractions of extremal rays of
a natural resolution of the nilpotent orbit closure.
\end{abstract}

\section{Introduction}
Let $\0$ be a nilpotent orbit in a simple complex  Lie algebra $\g$ (with $G$ its adjoint group).
The closure  $\overline{\0}$ in $\g$
 is a singular variety whose smooth part
admits a holomorphic symplectic Kostant-Kirillov form $\omega$.
A resolution $f: Z \to \overline{\0}$ is said {\em symplectic}
if $f^* \omega$ extends to a symplectic form on the whole of $Z$,
or equivalently if $f$ is crepant.
Note that there exist  nilpotent orbit closures which 
admit no symplectic resolutions(\cite{Fu1}).

A natural $G$-equivariant projective resolution of 
$\overline{\0}$  is given by $\mu: G \times^P \n \to
\overline{\0},$
where $P$ is a parabolic subgroup associated to the Jacobson-Morozov sub-algebra of
$\0$ and $\n$ is a nilpotent ideal in $\p = Lie(P)$ 
(see Section 2).
This resolution plays an important role in the study of 
singularities of $\overline{\0}$ (see for example \cite{Pan}).
 The resolution
$\mu$ is symplectic if and  only if $\0$ is an even orbit
(\cite{Fu3}). 
The closure $\overline{Amp}(\mu)$ of the ample cone of $\mu$ is
 a simplicial polyhedral cone and 
a face of $\overline{Amp}(\mu)$ corresponds
 to the Stein factorization of the natural birational map 
 $p: G \times^P \n \to G \times^Q (Q \cdot \n),$
for a suitable parabolic subgroup $Q$ of $G$ containing $P$. 
 Notice that $\mu $ factorizes through $p$, which gives
a birational map $\pi:  G \times^Q (Q \cdot \n) \to \overline{\0}$.
A natural question is if we can choose  $Q$ such that  $\pi$ is 
a symplectic resolution.  If it is the case and $\0$ is not an 
even orbit, 
then $p$ becomes an extremal contraction  of $\mu$,
 which is an important class of 
morphisms in Mori theory.  
We prove (Theorem \ref{classical},  Theorem \ref{exceptional},
Corollary \ref{exc} and Example \ref{D7}) 
that except for some particular orbits in $D_n$ and
 the orbit $\0_{D_7(a_2)}$ in  $E_8$,  
one can always obtain a symplectic
resolution of the nilpotent orbit closure $\overline{\0}$
 either by $\mu$ itself or by extremal contractions of $\mu$, 
provided that $\overline{\0}$ admits a symplectic resolution.
The proof is based on Proposition \ref{key} and some combinatorial calculations. Here and throughout the paper,
the notations of nilpotent orbits in exceptional Lie algebras are those in \cite{CM}.

Then we turn to study the birational geometry of Springer maps. 
Recall that (\cite{Ric}) for any parabolic sub-group $Q$ of 
$G$, the image of the moment map $T^*(G/Q)  \to \g \simeq \g^*$
is  a nilpotent orbit closure $\overline{\0}$. The orbit 
$\0$ will be called a {\em Richardson orbit} and $Q$ its 
{\em polarization}. The morphism $T^*(G/Q) \to \overline{\0}$
will be called the {\em Springer map} associated to $Q$, 
which is a generically finite  surjective projective map. 
For two Springer maps 
$T^*(G/Q_i) \to \overline{\0}, i=1, 2$ having the same degree, 
 we prove (Corollary \ref{bir}) that
there exists a birational map $ T^*(G/Q_1) \dasharrow T^*(G/Q_2)$ 
over $\overline{\0}$ which is decomposed into a finite sequence of 
(analytically) locally trivial families of stratified Mukai flops. 
Recall that (\cite{Fu1}) any symplectic resolution of 
$\overline{\0}$ is given by a (degree one) Springer map, so 
this implies that two symplectic resolutions of $\overline{\0}$
are connected by stratified Mukai flops, which has been previously
proved by Namikawa (\cite{Na}) using an ingenious argument.
 Our proof here is  different from \cite{Na}
 and does not make use of the
Springer correspondence and Mori theory, instead we will use a theorem of Hirai (\cite{Hi}) on polarizations of nilpotent orbits.
Here are three features of this result: 

(i) Nilpotent orbits usually have $G$-equivariant coverings, which turn to be important
in representation theory (e.g. \cite{BK}).
 An immediate corollary of our result is that Springer resolutions 
of coverings of nilpotent orbits are related by stratified Mukai flops. 

(ii) As a special case of Kawamata's conjecture that $K$-equivalence implies $D$-equivalence,
one conjectures the derived equivalence of the birational map  $T^*(G/Q_1) \dasharrow T^*(G/Q_2).$
As remarked in \cite{Na}, our result here reduces this conjecture to the cases of stratified Mukai flops. 

(iii) This result  provides evidence to the following: \vspace{0.2 cm}

{\bf Speculation:} Any two (projective) symplectic resolutions of a symplectic singularity
 are related by stratified Mukai flops. Any birational map between two projective
hyperK\"ahler manifolds is decomposed into a sequence of   stratified Mukai flops.
\vspace{0.2 cm}

Finally we turn to study stratified Mukai flops. Let 
$\phi: T^*(G/Q) \dasharrow T^*(G/Q')$ be such a flop. We show (Theorem \ref{Mukai}) that the graph closure
of $\phi$ is isomorphic to the variety $G \times^P \n$ which gives the natural 
resolution of $\overline{\0}$ and the two graph projections $ T^*(G/Q) \leftarrow G \times^P \n \to T^*(G/Q')$
are given by the two 
contractions of extremal rays  of $\mu$. This gives a conceptional and explicit description of the flop $\phi$.

Some interesting by-products are obtained by our methods here.
We prove (Theorem \ref{exceptional} and Corollary \ref{exc}) that for 
the orbits $D_4(a_1)+A_1$ in $E_7$, $D_6(a_1)$ and $D_7(a_2)$ in $E_8$,
 their closures admit a symplectic resolution, while for the four orbits $A_4+A_1, D_5(a_1)$ in $E_7$,  
$E_6(a_1)+A_1$ and $E_7(a_3)$ in $E_8$, their closures do not admit a
symplectic resolution (Corollary \ref{exc}). Together with \cite{Fu1} (see also 
Proposition \ref{class_f}, Proposition \ref{exc_f}),
this completes the classification  
of nilpotent orbits in a simple Lie algebra whose closure admits a symplectic resolution.

As remarked in \cite{BM}, the degree of the Springer map 
$ T^*(G/Q) \to \overline{\0}$ plays an important role in several different contexts 
(e.g. \cite{BB}, Theorem 5.5, 5.6, 5.8).
Another by-product of this paper is that we can determine (Remark \ref{det}) the degree of the Springer map
associated to any parabolic sub-group $Q$
in a very explicit and practical way.
When $\g$ is classical, this is due to Hesselink (\cite{He}). When $\g$ is exceptional, this result seems to
be new.  
 \vspace{0.2 cm}

{\em Notations:} We fix a Cartan sub-algebra $\h$, a Borel
 sub-algebra $\b$, a system of positive roots $\Phi^+$ and
 simple  roots
$\Delta = \{ \alpha_1, \cdots, \alpha_n \}$. The labels of roots
in the Dynkin diagram 
are the same as those in \cite{Bou}.  For a subset $\Gamma \subset \Delta$, we denote by
$\p_\Gamma$  the standard parabolic sub-algebra 
$\mathfrak{b} \oplus_{\beta \in \langle \Delta - \Gamma \rangle^- }
 \g_\beta$ and $P_\Gamma$ the standard parabolic subgroup in $G$
with Lie algebra $P_\Gamma$.
Note that in the literature, our $\p_\Gamma$ is usually denoted by $\p_{\Delta - \Gamma}.$
The marked Dynkin diagram of $P_\Gamma$ is obtained by marking the nodes in $\Gamma$. We will denote by $\u(P_\Gamma)$ or $\u(\p_\Gamma)$
the nil-radical of $\p_\Gamma$ and by $\l(P_\Gamma)$ 
or $\l(\p_\Gamma)$ the Levi factor of $\p_\Gamma$.
More precisely, $\u(\p_\Gamma) =
 \oplus_{\Phi^+ -  \langle \Delta - \Gamma \rangle^+} \g_\beta$ 
and  $\l(\p_\Gamma) = 
\h \oplus_{ \beta \in \langle \Delta - \Gamma \rangle } \g_\beta.$
\vspace{0.2 cm}

{\em Acknowledgements:} Part of this work has been done during my visit to I.M.S. 
at the Chinese University of Hong Kong and  M.S.R.I. at Berkeley.
 It's my pleasure to thank C. Leung and Y. Ruan for the kind invitations and the two
institutes for their hospitality. I want to thank L. Manivel for  the references \cite{Dem} \cite{NSZ},
and D. Alvis for explanations on \cite{Alv}.

\section{Extremal contractions}

Let $\g$ be a simple  Lie algebra and $G$ its adjoint group. 
For a nilpotent element $x \in \g$, the Jacobson-Morozov theorem
gives an $\mathfrak{sl}_2$-triplet $(x, y, h)$, i.e. 
$[h, x] = 2 x, [h, y] = -2y, [x, y] = h$. Up to replacing this
triplet by a conjugate one, we can assume that $h \in \h$ and $h$ 
is $\Delta$-dominant.
This triplet makes $\g$ an $\mathfrak{sl}_2$-module, so we have a 
decomposition $\g = \oplus_{i \in \zit} \g_i$, where
 $\g_i = \{z \in \g \mid [h, z] = i z \}.$ The Jacoboson-Morozov 
parabolic sub-algebra of this triplet is 
$\p: = \oplus_{i \geq 0} \g_i$. Its conjugacy class is uniquely 
determined by the  the nilpotent orbit $\0 = G \cdot x$.

Recall that a nilpotent orbit is uniquely determined  by its weighted Dynkin 
diagram, which is obtained by assigning $\alpha(h)$ to the node $\alpha$. Under
our assumption, $\alpha(h) \in \{0, 1, 2\}$.
\begin{Prop}\label{JM}
The marked Dynkin diagram of $\p$ is obtained from the weighted 
Dynkin diagram of $\0$ by marking the nodes with non-zero weights.
\end{Prop}
\begin{proof}
Let $\Gamma$ be the set of marked nodes in the marked Dynkin diagram
of $\p$, then the set 
 $\Delta - \Gamma$ consists of simple roots $\alpha$ such that
$\g_\alpha, \g_{- \alpha}$ are contained in $\p$. Notice that
$\g_\alpha \subset \g_{\alpha(h)}, 
\g_{- \alpha} \subset \g_{- \alpha(h)}$, 
which gives that $\alpha(h) = 0.$
\end{proof}

The closure $\overline{\0}$ of $\0$ in $\g$ is singular. A natural 
resolution of $\overline{\0}$ is given by 
$\mu: G \times^P \n  \to \overline{\0}$, where $P$ is a connected
subgroup of $G$ with Lie algebra $\p$ and 
$\n := \oplus_{i \geq 2} \g_i$ is a nilpotent ideal of $\p$.
Note that the variety $ G \times^P  \n$ and the map $\mu$
is independent (up to isomorphisms)
 of the choice of the element $ x \in \0$ and the 
standard $\mathfrak{sl}_2$-triplet.
This resolution is symplectic
if and only if $\0$ is an {\em even orbit}, i. e.  $\g_1 = 0$, or
equivalently the weights in the weighted
Dynkin diagram of $\0$  are only $0$ and $ 2$ (\cite{Fu3}). 

 If we denote by 
$\tilde{\mu}$ the Stein factorization of $\mu$, then 
$\overline{Amp}({\tilde{\mu}})$ and
 $\overline{NE}({\tilde{\mu}})$ are both 
 simplicial polyhedral cones. 
The contraction of a face in 
$\overline{NE}({\tilde{\mu}})$ is given by the Stein
 factorization of the morphism
$$ p:  G \times^P \n  \to G \times^Q (Q \cdot \n),$$
for a suitable  parabolic
sub-group $Q$ in $G$ containing $P$.
The map $\mu$ factorizes through $p$, which gives a birational
map  $\pi:  G \times^Q (Q \cdot \n)  \rightarrow  \overline{\0}$.
 An interesting  question is if we can  choose
$Q$ such that $\pi$ becomes a symplectic resolution. 
\begin{Rque} \upshape
When $\0$ is even, then $\mu$ is already a symplectic resolution.
In this case, any  contraction as above will produce
a singular variety $ G \times^Q (Q \cdot \n)$. 
\end{Rque}

\begin{Prop} \label{key} 
Let $Q$ be a parabolic subgroup containing $P$ and
 $\u(Q)$  the nil-radical of $\q: = Lie(Q)$.  
If $\n \subset \u(Q)$ and $ 2 \dim (\u(Q)) = \dim \0$,
then 

(1) $ G \times^Q (Q \cdot \n) \simeq T^*(G/Q)$ and $\pi$
is a symplectic resolution,

(2) if $\0$ is not an even orbit, then the map $p$
is a contraction of an extremal face of $\overline{NE}({\tilde{\mu}})$.
\end{Prop}
\begin{proof}
By assumption, we have $Q \cdot \n \subset Q \cdot \u(Q) = \u(Q)$.
Note that $Q \cdot \n$ is the image of $Q \times^P \n$
under the projective map $Q/P \times \g \to \g$, so it is closed
in $\g$. Since $\pi$ is birational, we have 
$\dim (Q \cdot \n) = \dim (\0) - \dim(G/Q) = \dim(\0) - \dim(\u(Q))
= \dim(\u(Q))$, which gives $Q \cdot \n = \u(Q)$ since $\u(Q)$
is closed and irreducible.  
Now assertion (1) follows immediately.

Assume that $\0$ is not even, then $p$ is not an isomorphism. The exceptional set 
$E$ of $p$ has pure codimension 1 since $ G \times^Q (Q \cdot \n) $ is smooth.
Let $E = \cup_i E_i$ be the decomposition into irreducible components, then 
$K:= K_{G \times^P \n }  = \sum_i a_i E_i$ with $a_i >0$. Now it is easy to see that
$K \cdot C < 0$ for any curve $C$ contracted to a point by $p$. For any such a curve $C$, its
class lies in an extremal face of $\overline{NE}({\tilde{\mu}})$, the one dual to the face in
$\overline{Amp}(\tilde{\mu})$ determined by the map  $p$.
\end{proof}
\begin{Rque}\upshape
The advantage of using extremal contractions is that
the degree of the map $T^*(G/Q) \to \overline{\0}$ is automatically 1.
 This turns out to be helpful when $\g$ is exceptional,
since  in this case, it is not easy to calculate the degrees of Springer
maps.
\end{Rque}

One should bear in mind that even when $T^*(G/Q) \to \overline{\0}$
is a symplectic resolution, in general the birational map 
$ G \times^Q (Q \cdot \n)  \to \overline{\0}  $ is not a symplectic 
resolution. In fact, the variety $ G \times^Q (Q \cdot \n) $
can  even be singular if we drop the conditions in the precedent proposition,
 as shown by the following example.

\begin{Exam}\upshape
In $\mathfrak{sl}_5$, let $\0$ be the nilpotent orbit with Jordan 
type $[4, 1]$. Then the Jacobson-Morozov parabolic subgroup  of $\0$
is a Borel subgroup. Let $Q$ be the standard parabolic subgroup with
flag type $[2, 1, 1, 1]$, then $T^*(G/Q) \to \overline{\0}$ is
a symplectic resolution, so $\dim \u(Q) = 1/2 \dim \0$, but 
$\n $ is not contained in $\u(Q)$. A direct calculus shows that
$Q \cdot \n$ is defined by some quadric equations, and it is singular in
codimension 1.
\end{Exam}

\section{ Classical types}

For a nilpotent orbit $\0$ in a classical 
simple Lie algebra $\g$, we denote by
$\Theta_i$ the set of  nodes in $\Delta$ with weight $i$, 
for $i = 0, 1$ or $2$.  The standard Jacobson-Morozov
parabolic sub-algebra of $\0$ will be denoted by $\p$, which is
obtained by marking $\Theta_1 \cup \Theta_2$ (Prop. \ref{JM}).
We denote by $P$ the connected subgroup of $G$ with Lie algebra
$\p$. Let ${\bf d} = [d_1, \cdots, d_k]$ be the Jordan type of 
$\0$. Recall the following classification theorem from \cite{Fu1}.
\begin{Prop}\label{class_f}
Assume $\g$ is simple classical. Then 
the closure $\overline{\0}$ admits a symplectic resolution if and only
if $\0$ is in the following list:

(i) $\0$ is in $\g = \mathfrak{sl}_n$;

(ii) $\g = \mathfrak{so}_{2n+1}$ (resp. $\mathfrak{sp}_{2n}$) and
there exists an odd (resp. even) number $q \geq 0$ such that
$d_1, \cdots, d_q$ are odd and $d_{q+1}, \cdots,d_k$
are even;  

(iii) $\g = \mathfrak{so}_{2n}$ and  either 

(iii-a) there exists an even
number $q \geq 4$ such that $d_1, \cdots, d_q$ are odd and 
$d_{q+1}, \cdots,d_k$ are even;  or 

(iii-b) there exist exactly two 
odd parts in ${\bf d}$ at positions $2t-1$ and $2t$ for some number
 $t \geq 1$. 
\end{Prop}

From now on, we will assume that $\0$ is not an even orbit, i. e. the Jordan
type of $\0$ has parts with different parities.
We  have  the following lemma, whose proof is left to the reader.
\begin{Lem}
For cases (i), (ii) and (iii-a) in Proposition \ref{class_f}, 
the set  $\Theta_1$ has even number of elements.
For case (iii-b) in Proposition \ref{class_f}, 
there are the following cases:

(iii-b-1) if $t = 1$, then $\alpha_n \in \Theta_1$ and the number of elements in 
$\Theta_1$ is odd;

(iii-b-2) if $k = 2t \geq 4$, then the number of elements in 
$\Theta_1$ is even and $\alpha_{n-1}, \alpha_n \in \Theta_1$;

(iii-b-3) if $k > 2t \geq 4$, then the number of elements in $\Theta_1$ is even.
\end{Lem}

For cases (i), (ii) and (iii-a), 
we decompose the set
$\Theta_1  = :\{ \alpha_{m_i}| i = 1, \cdots, 2l, m_i < m_{i+1}, \forall i \}$
as the disjoint union $\Theta_1^I \cup \Theta_1^{II}$, where
$\Theta_1^I = \{ \alpha_{m_{2j-1}}| j = 1, \cdots, l \}  $
 and $\Theta_1^{II} = \{ \alpha_{m_{2j}}| j = 1, \cdots, l \}.  $
For case  (iii-b-1), we obtain a decomposition in a similar way 
 $\Theta_1  = \Theta_1^I \cup \Theta_1^{II}$, with the extra
element 
$\alpha_n$ in $\Theta_1^{II}$.
For case  (iii-b-2), we have two distinct 
 decompositions (except when $d_{2t-1} = d_{2t} = 1$) of  $\Theta_1$
 as   $\Theta_1^I \cup \Theta_1^{II} = 
 \Theta_1^{' I} \cup \Theta_1^{' II}, $
where $\Theta_1^I, \Theta_1^{II}$ are defined as above and
$\Theta_1^{' I}$ (resp.  $\Theta_1^{' II}$  ) is obtained from  $ \Theta_1^I$
by replacing $\alpha_{n-1}$ (resp. $\alpha_n$) by  $\alpha_n$ (resp. $\alpha_{n-1}$).

Let $\q_1$ (resp. $\q_2$, $\q'_1$, $\q'_2$ )
be the standard parabolic sub-algebra obtained by marking the nodes
 in
$\Theta_1^I \cup \Theta_2$ (resp. $\Theta_1^{II} \cup \Theta_2$, 
$\Theta_1^{' I} \cup \Theta_2$, $\Theta_1^{' II} \cup \Theta_2$).
Let $Q_1, Q_2, Q'_1, Q'_2$ be the parabolic subgroups with Lie algebras
$\q_1, \q_2, \q'_1, \q'_2$ respectively.
One remarks that the standard Jacobson-Morozov parabolic sub-algebra
$\p$ is contained in theses sub-algebras.
 Let 
$\pi_i: G \times^{Q_i} (Q_i \cdot \n) \to \overline{\0}, i=1, 2,$ 
be two maps obtained from contractions of $\mu$. Similarly 
one has $\pi'_i$ for case (iii-b-2).  

\begin{Thm} \label{classical}

For case (i), the two maps $\pi_1, \pi_2$ are both symplectic resolutions.
 The rational map $\pi_2^{-1} \circ \pi_1 $ is resolved by $G \times^P \n$;

For cases (ii), (iii-a) and (iii-b-1) 
the map $\pi_2$ is a symplectic resolution;

For case (iii-b-2), $\pi_i, \pi'_i$ are all symplectic resolutions, which are
all dominated by $\mu$.

For case (iii-b-3), none of the maps arising
from extremal contractions of $\mu$ is a symplectic resolution. 
\end{Thm}
\begin{proof}

The strategy is to apply Prop. \ref{key}, so we need to check the conditions
$\n \subset \u(\q)$ and  $ \dim (\u(\q)) = 1/2 \dim \0$. The first condition
is easily checked by our choice of the decomposition of $\Theta_1$. 
The second condition is equivalent to
 $\dim \u(\p) - \dim  \u(\q) = \dim  \u(\q) - \dim \n$. We will check this
condition case by case.

For case (i), 
$\dim \u(\p) - \dim  \u(\q_1)$ is the number of 
positive roots $\beta = \sum_{i \leq k \leq j} \alpha_k$
such that there exists a unique $k_0$ such that 
$\alpha_{k_0} \in \Theta_1^{II} $ and for $i \leq k \neq  k_0 \leq j$, 
we have  $\alpha_k \in \Theta_0$. This is also the sum
$\sum_{k=1}^l N(\alpha_{m_{2k}})$, where $N(\alpha_{m_{2k}})$
is the number of connected subgraphs containing the node 
$\alpha_{m_{2k}}$ and the other nodes are in $\Theta_0$.
On the other hand,
$\dim  \u(\q_1) - \dim \n $ 
is the number of positive roots 
$\beta = \sum_{i \leq k \leq j} \alpha_k$
such that the sum of weights of nodes in $\beta$ is 1 and 
there exists some $i \leq k_0 \leq j$ such that $\alpha_{k_0}
\in \Theta_1^I$. This number is the sum 
$ \sum_{k =1}^l N(\alpha_{m_{2k-1}}) $ . Note that 
the weighted Dynkin diagram of any nilpotent orbit in
 $\mathfrak{sl}_n$ is invariant under the non-trivial
graph automorphism (Lemma 3.6.5 \cite{CM}), so we have
$N(\alpha_{m_{2k-1}}) = N(\alpha_{m_{2l+2-2k}})$, which gives 
the equality 
$\dim \u(\p) - \dim  \u(\q_1) = \dim  \u(\q_1) - \dim \n$. Similar arguments 
apply to $\q_2$. Thus $\pi_1, \pi_2$ are both symplectic and dominated by $\mu$.

For case (ii), first consider  $\g = \mathfrak{so}_{2n+1}$.
As easily seen, the weighted Dynkin diagram of $\0$ has the
 following form (where nodes are replaced by their weights):
$$ \cdots- 2 - 0^{2a} - 1 - 0^{r_1} - 1 -  0^{2a} - 1 - \cdots - 
0^{2a} -1 - 0^{r_l} - 1 - 0^{a-1} \Rightarrow 0 $$
where $r_1, \cdots, r_l$ are non-negative integers, 
$2a+1 = q$ and $0^m$ means the consequentive $m$
nodes have weights $0$. Note that the weights of nodes on the left-hand side of
the node with weight 2 can be only 0 or 2.
Like  case (i), one has
$$\dim \u(\p) - \dim  \u(\q_2) = \sum_{i=1}^l N(\alpha_{m_{2i-1}})
= \sum_{i=1}^l (2a+1) (r_i +1) = q (\sum_{i=1}^l r_i + l).$$
On the other hand 
$$ \dim  \u(\q_2)  - \dim  \n =   
\sum_{i=1}^{l-1} N(\alpha_{m_{2i}}) + N = 
 q (\sum_{i=1}^{l} r_i + l) + (N - q (r_l + 1)),$$
where $N$ is the number of positive roots
such that one (with multiplicity)  of whose summands is $\alpha_{n-a}$ (the rightmost
node with weight 1)
 and the others
are in $\{\alpha_j | n-a-r_l  \leq j \neq n-a  \}$.
Using the table for positive roots (\cite{Bou}), we
find $N = q (r_l + 1)$, which gives the assertion.

The proof for $\g = \mathfrak{sp}_{2n}$ is similar to the case of 
$\mathfrak{so}_{2n+1}$.
The key point is to notice that the weighted Dynkin diagram of 
$\0$ is of the following form:
$$ \cdots -2 - 0^{q-1} - 1 - 0^{r_1} - \cdots - 0^{q-1} - 1 - 
0^{r_l} - 1 - 0^{q/2 -1} \Leftarrow 0,$$ where $r_1, \cdots, r_l$
are non-negative integers.
For case (iii-a), the argument is the same, by noticing that the weighted
Dynkin diagram has the following form: 
$$ \cdots - 2 - 0^{q-1} -1 - 0^ {r_1} -1 - \cdots - 0^{r_l} - 0^{q/2},$$
where $0^{q/2}$ means the last $q/2$ nodes have weights 0. 
The other two cases (iii-b-1) and (iii-b-2) are similar.

For case (iii-b-3), we notice that the Levi type
of any (degree 1) polarization $Q$  of $\0$ is
 ${\bf d'}= [d_1, \cdots, d_{2t-2}, d_{2t-1}+1, d_{2t}-1, d_{2t+1},
 \cdots, d_k].$  The dual partition of  ${\bf d'}$ has the form
$ord({\bf d'}) := [k^{2s},q_1, \cdots, q_l]$, with 
$q_j \leq k-2s, \forall j$ and $s \geq 1$, since
even parts appear with even multiplicity in the Jordan type of any nilpotent
orbit in $\mathfrak{so}_{2n}$. By our assumption, $q_l \geq 2$ since $d_1$ is even. 
By \cite{He}, this implies that every flag type determined by $ord({\bf d'})$
corresponds to two marked Dynkin diagrams, i.e.  either $\alpha_{n-1}$ or $\alpha_n$ is marked, but not
both. 

 The weighted Dynkin diagram of $\0$ has the following form:
$$ \cdots - 1 - 0 - 1 - 0^{k-3} - 1^2,$$
where $1^2$ means that the nodes $\alpha_{n-1}, \alpha_n$
have weights 1.  Suppose that $Q$  
contains $P$, then in the flag type of $Q$, there is either a part equal to $k-1$ or a part
equal to $q' \geq k+1$, which is a contradiction. 
\end{proof}

\begin{Rque}
(1). Notice that for $\mathfrak{so}_{2n+1}$,  although one  has
$\n \subset \u(\q_1)$, but 
$\dim \u(\p) - \dim  \u(\q_1) > \dim  \u(\q_1) - \dim \n$, 
thus $\pi_1$ is not a symplectic resolution. Similar remark applies to other cases.

(2).  For case (iii-b-3), as we know the Levi type 
of a degree one polarization, we can also obtain the marked Dynkin diagram of a symplectic resolution. Later on, 
we will find all symplectic resolutions of a nilpotent orbit closure by starting from any given one.  
\end{Rque}

\section{Exceptional cases}

Let us recall the following classification result from \cite{Fu1}.
The notations of orbits are those in \cite{CM} (p. 128-134).

\begin{Prop}\label{exc_f}
(i) For the following Richardson orbits, we do not know if 
their closures admit a symplectic resolution or not:

 $D_4(a_1)+A_1, A_4+A_1, D_5(a_1)$  in $E_7$ and

 $D_6(a_1), D_7(a_2), E_6(a_1)+A_1, E_7(a_3)$ in $E_8$.

(ii) For other orbits in a simple exceptional Lie algebra $\g$,
its closure admits a symplectic resolution if and only if it is a 
Richardson orbit. The following is the complete list of 
such orbits:

(ii-a) even orbits;

(ii-b) $C_3$ in $F_4$,  $2A_1, A_2+2A_1, A_3, A_4+A_1, D_5(a_1)$
in $E_6$, $D_5+A_1, D_6(a_1)$ in $E_7$ and $A_4+A_2+A_1, A_6+A_1,
E_7(a_1)$ in $E_8$.
\end{Prop}

\begin{Thm} \label{exceptional}
For the orbits in (ii-b) of Proposition \ref{exc_f}
 and orbits $D_4(a_1)+A_1$ in $E_7$,
$D_6(a_1)$ in $E_8$, one can always obtain a
symplectic resolution by an extremal contraction of the natural 
resolution  $\mu: G \times^P \n \to \overline{\0}$.
In particular, the closures of the orbit  $D_4(a_1)+A_1$ in $E_7$
and the orbit
$D_6(a_1)$ in $E_8$ admit a symplectic resolution.
\end{Thm}
\begin{proof}
We will verify the conditions in Proposition \ref{key}, and then 
apply it to conclude.  The notations of roots are those 
in \cite{Bou}. We just give the corresponding polarization $Q$
to each orbit in the statement.  The condition
$\n \subset \u(\q)$ and  the dimension check can be done
by using the tables of root systems in \cite{Bou}. 
For example for $\g = F_4$ and $\0 = \0_{C_3}$, the Jacobson-Morozov
standard parabolic is $P = P_{\alpha_1 \alpha_2 \alpha_3}$ while $Q =  P_{\alpha_3 \alpha_4}.$
One checks that $\dim \u(P) - \dim \u(Q) = \dim \u(Q) - \dim \n = 2.$
Notice that if we take $Q' = P_{\alpha_1 \alpha_4}$, then $\n$ is not contained
in $\u(Q')$, since $\alpha_2 + 2 \alpha_3$ is still a positive root, so Proposition \ref{key} is not applicable to $Q'$. In the tables
below, the first row gives the Lie algebras, the second row lists
the nilpotent orbits and the third row gives the corresponding (degree one)
polarizations.

\begin{tabular}{|l|l|l|l|l|l|}
\hline 
$F_4$  & \multicolumn{5}{l|}{$E_6$}  \\
\hline
$C_3$ & $2A_1$  & $A_2+2A_1$   & $A_3$  & $A_4+A_1$ & $D_5(a_1)$ \\
\hline
$ P_{\alpha_3 \alpha_4}$ & $P_{\alpha_1}, P_{\alpha_6}$
&  $P_{\alpha_3}, P_{\alpha_5}$  & $P_{\alpha_1 \alpha_2},
P_{\alpha_2 \alpha_6} $ &  $P_{\alpha_3 \alpha_5}$ & 
 $P_{\alpha_2 \alpha_3 \alpha_6}, 
P_{\alpha_1 \alpha_2 \alpha_5} $ \\
\hline
\end{tabular} \vspace{0.2 cm}

\begin{tabular}{|l|l|l|l|l|l|l|}
\hline 
 \multicolumn{3}{|l|}{$E_7$}  & \multicolumn{3}{l|}{$E_8$}  \\
\hline
$D_4(a_1)+A_1$  & $D_5+A_1$ &  $A_4+A_2+A_1$ & $A_6+A_1$ &
$ A_6+A_1$  &  $D_6(a_1)$   \\
\hline
$ P_{\alpha_2 \alpha_7}$ & $P_{\alpha_1 \alpha_3 \alpha_5}$
& $P_{\alpha_1 \alpha_2 \alpha_3 \alpha_7}$ & $P_{\alpha_3}$
&$P_{\alpha_4}$ & $P_{\alpha_1 \alpha_2 \alpha_3}$\\
\hline
\end{tabular} \vspace{0.2 cm}

Finally for the orbit 
$ E_7(a_1)$ in $E_8$, we take
 $Q = P_{\alpha_1 \alpha_2 \alpha_3  \alpha_7 \alpha_8}.$

\end{proof}

\begin{Rque}
The proof also gives another way to show that these
orbits are Richardson. 
\end{Rque}

\section{Birational geometry}

The precedent sections give a particular symplectic resolution
of a nilpotent orbit closure $\overline{\0}$ 
provided we know the existence of such a resolution.
In this section, we will describe
a way to find all symplectic resolutions of $\overline{\0}$ starting from any given one. This procedure has been previously
described in \cite{Na}.

For two standard parabolic subgroups $P_\Gamma$ and $P_{\Gamma'}$, 
we define $P_\Gamma \sim_R P_{\Gamma'}$ 
(or $\Delta - \Gamma \sim_R  \Delta - \Gamma'$)
 if the Richardson orbits
corresponding to  $P_\Gamma$ and  $P_{\Gamma'}$ are the same, say 
$\0$. We say  that $P_\Gamma$ and  $P_{\Gamma'}$ are {\em equivalent}
(write $P_\Gamma \sim P_{\Gamma'}$ or  
$\Delta - \Gamma \sim \Delta - \Gamma'$ )
if furthermore the degrees of the two Springer
 maps $T^*(G/ P_\Gamma) \to \overline{\0}
\leftarrow T^*(G/ P_{\Gamma'})$ are the same. 

\begin{Thm}[Hirai \cite{Hi}] \label{Hi}
Assume $\g$ is simple.
The equivalence relation $\sim_R$ 
is generated by the following fundamental ones:

(1) In $B_n$ or $C_n$ with $n=3k - 1, k\geq 1$, 
$P_{\alpha_{2k-1}} \sim_R P_{\alpha_{2k}}.$

(2) In $D_4$, 
$ P_{\alpha_2} \sim_R P_{\alpha_{3} \alpha_4}.$

(3) In $D_n$ with $n = 3k +1, k\geq 2$, 
$P_{\alpha_{2k}} \sim_R P_{\alpha_{2k+1}}.$

(4) In $G_2$, $P_{\alpha_1} \sim_R P_{\alpha_2}.$

(5) In $F_4$, $P_{\alpha_2} \sim_R P_{\alpha_3}  
\sim_R P_{\alpha_1 \alpha_4}  $.

(6) In $E_6$, $P_{\alpha_4}  \sim_R P_{\alpha_2 \alpha_5}  $.
 
(7) In $E_8$, $P_{\alpha_5}  \sim_R P_{\alpha_2 \alpha_3}  $.

(8) In $A_n$,  $P_{\alpha_i} \sim_R P_{\alpha_{n+1-i}}, \forall i$.

(9) In $D_{2k+1} (k\geq 2)$, 
$P_{\alpha_{2k}} \sim_R P_{\alpha_{2k+1}}.$

(10) In $E_6$, $P_{\alpha_1} \sim_R P_{\alpha_6}$,
and $P_{\alpha_3} \sim_R P_{\alpha_5}$.

(GP)[General principle] If $\Delta_1, \Delta_2$ 
are two subsets of  $\Delta$ orthogonal to each other.
Let $\Gamma_i \subset \Delta_i, i=1, 2$ be two subsets
and $\Gamma'_1 \subset \Delta_1$ a subset such that
$\Delta_1 - \Gamma_1 \sim_R \Delta_1 - \Gamma'_1 $
in the root system $\langle \Delta_1 \rangle$, then 
$P_{\Gamma_1 \cup \Gamma_2} \sim_R P_{\Gamma'_1 \cup \Gamma_2}    $. Here $\Delta_i$ can be empty.
\end{Thm}

The proof of this theorem is essentially a type-by-type check, since one can 
determine the Richardson orbit of any parabolic subgroup (in classical cases,
this is given by the Spaltenstein map, while in most exceptional cases, it
suffices to do just a dimension counting. Some particular attention should be payed
to a few cases, for details see \cite{Hi}).

Our result is to give a list of fundamental relations for
the equivalence $\sim$.
\begin{Thm}\label{main}
Suppose that $\g$ is simple. Then the equivalence  $\sim$ 
is generated by the relations (8), (9), (10) and (GP) in Theorem
\ref{Hi}.
\end{Thm}
\begin{Rque} \upshape
As we can see in the proof, the theorem is not true
if $\g$ is not simple. 
\end{Rque}

We begin the proof  by some lemmas.
\begin{Lem} \label{equiv}
(i) For the parabolic subgroups appeared 
in each equivalence relation of  (1) - (7) in Theorem \ref{Hi},
 there is only
one such that the associated Springer map is birational. 

(ii) For any parabolic subgroup appeared
in the relations (8) - (10) in Theorem \ref{Hi},
the associated Springer map is birational.
\end{Lem}
\begin{proof}
For the case of $B_{3k-1}$ (resp. $C_{3k-1}$), 
the Richardson orbit has
Jordan type $[3^{2k-1}, 1^2]$ (resp. $[3^{2k-2}, 2^2]$).
 The Springer map
associated to $P_{\alpha_{2k-1}}$ (resp. $P_{\alpha_{2k}}$) 
 is of degree 1, while that
of  $P_{\alpha_{2k}}$ (resp.  $P_{\alpha_{2k-1}}$) is of degree 2. 
Here we used Hesselink's formula for the degrees of 
Springer maps (Theorem 7.1 \cite{He}) in classical Lie algebras.

In $D_4$, the Richardson orbit of $P_{\alpha_2}$
and   $P_{\alpha_{3} \alpha_4}$
 has Jordan type $[3^2, 1^2]$. One calculates
that the degrees of the Springer maps
associated to  the two polarizations are respectively
1 and 2. 

In $D_{3k +1}, k\geq 2$, the Richardson orbit for 
$P_{\alpha_{2k}}$ and  $P_{\alpha_{2k+1}}$ has Jordan type
$[3^{2k}, 1^2]$. The degrees of the
Springer maps are respectively 1 and 2. 

In $G_2$, the Richardson orbit of  $P_{\alpha_1}$ 
and $P_{\alpha_2}$ is the sub-regular orbit $\0 := \0_{G_2(a_1)}.$
It is an even orbit with weighted Dynkin diagram
$ 2 \equiv > 0$, so $\pi_2$  is a symplectic 
resolution, where 
$\pi_i: T^*(G/P_{\alpha_i}) \to \overline{\0}, i=1, 2$ are the 
Springer maps. The closure $\overline{\0}$ is normal since
$\0$ is the sub-regular orbit, whose singular part
is the closure 
of the codimension 2 orbit $\0' := \0_{\tilde{A}_1}$. A slice
transversal  to $\0'$ has an isolated  normal surface singularity, 
which admits a unique crepant resolution.
Suppose that $\pi_1$ is  birational, then the birational
map $\phi: = \pi_2^{-1}  \circ \pi_1$ is an isomorphism 
over the pre-images of  $\0'$. 
 Let $\L$ be a $\pi_2$-ample line bundle and $C'$ an irreducible
component  of the $\pi_1$-fiber of a point in $\0'$, then the 
line bundle $\phi^*(\L)$ satisfies 
$\phi^*(\L) \cdot C' = \L \cdot \phi(C') > 0.$ 
But the Picard group of $T^*(G/P_{\alpha_1})$ is $\zit$ and $\pi_1$ is projective,
so $\phi^*(\L)$ is $\pi_1$-ample, which implies that $\phi$
is in fact an isomorphism. Note that $\pi_1, \pi_2$ are
both $G$-equivariant, so is the isomorphism $\phi$.
This implies that $G/P_{\alpha_1}$ and $G/P_{\alpha_2}$ 
are isomorphic as $G$-varieties, 
which is absurd since $P_{\alpha_1}$ and $P_{\alpha_2}$
are not $G$-conjugate. 
In conclusion,  $\pi_1$ is not birational. 

Now consider (5). The Richardson orbit is given by
$\0 = \0_{F_4(a_3)}$, whose weighted Dynkin diagram is given by
$0 - 2 \Rightarrow 0 - 0.$ This is an even orbit, thus
the Jacobson-Morozov parabolic subgroup $P_{\alpha_2}$
gives a symplectic resolution.  
Notice that the Picard group of $G/P_{\alpha_1 \alpha_4}$ is
$\zit^2$ which is different from that of $G/P_{\alpha_2}$,
so the Springer map associated to $P_{\alpha_1 \alpha_4}$ 
is not birational.
By \cite{Br}, the orbit closure $\overline{\0}$ is normal, 
whose singular part contains a codimension 2 orbit $C_3(a_1)$.
Now a similar argument as that for the relation (4) shows that
the Springer map of $P_{\alpha_3}$ is not birational.

For case (6), the Richardson orbit is $D_4(a_1)$, which is
an even orbit. The Jacobson-Morozov parabolic sub-group
 is $P_{\alpha_4}$,
so it gives a symplectic resolution, while 
$P_{\alpha_2 \alpha_5}$ does not, for the reason of different
Picard groups.

For case (7), the Richardson orbit is $E_8(a_7)$, 
which is again an even orbit with  the Jacobson-Morozov parabolic
sub-group
$P_{\alpha_5}$. The situation is similar to (6).

For case (8), (9), the Springer map associated to
each parabolic sub-group is birational, by Theorem \ref{classical}
(cases (i) and (iii-b-2)). For case (10), the associated Springer
map is birational by the proof of Theorem \ref{exceptional}.
\end{proof}
\begin{Rque} \upshape
(i) In the appendix, we will calculate explicitly the degree 
of the Springer map associated to each parabolic sub-group
appeared in $(4)-(7)$ of Theorem \ref{Hi}
 by using a formula of Borho-MacPherson
(\cite{BM}). However, we prefer to give the more geometric
 proof here.

(ii)  The two varieties $G_2/P_{\alpha_1}$ and  $G_2/P_{\alpha_2}$
are not isomorphic even as algebraic varieties, since they have 
different automorphism groups (see \cite{Dem}, also \cite{BK}).
The variety $F_4/P_{\alpha_2}$ is not isomorphic 
to $F_4/P_{\alpha_3}$ as algebraic varieties
since their Chow groups are different (see for example \cite{NSZ}).
\end{Rque}

Let $P \subset Q$ be two standard parabolic subgroups in $G$ with Lie
algebras $\p, \q$ and $L$ a Levi subgroup of $Q$. 
The projection to the first factor of the
direct sum $\u(\p) = \u(\l(\q) \cap \p) \oplus \u(\q)$
gives an $L$-equivariant map 
$f: L \cdot \u(\p) \to L \cdot \u(\l(\q) \cap \p)$.
Let $g: L \times^{(L \cap P)}  \u(\l(\q) \cap \p) \to  
L \cdot \u(\l(\q) \cap \p)$ and 
$g': Q \times^P   \u(\p)  \to Q \cdot  \u(\p) = L \cdot \u(\p)$
be the natural morphisms. Note that $g$ is a product
of  isomorphisms
with Springer maps (in some simple  Lie sub-algebras of $\l(\q)$),
 so it is generically finite. 
\begin{Lem} \label{lemdeg}
The morphism $g'$ is the pull-back via $f$ of the map $g$, i.e. the following
diagram is Cartesian. In particular, $\deg(g') = \deg (g)$.
$$\begin{CD}
Q \times^P \u(\p) @>g'>> Q \cdot \u(\p)=L\cdot \u(\p) \\
@VVV   @VfVV \\
L \times^{L \cap P} \u(\l(\q) \cap \p)  @>g>> L \cdot \u(\l(\q) \cap \p). 
\end{CD}$$

\end{Lem}
\begin{proof}
Let $Z$ be the fiber product of $f$ and $g$, then
we have an $L$-equivariant map 
$\eta: Z \to L/(L \cap P)$. The fiber of $\eta$ over the identity 
is isomorphic to $\u(\p)$.
If we denote by $U$ the unipotent subgroup
of $Q$, then $U \subset P$ and $Q/P = LU/P \simeq L/(L \cap P)$.
This shows that $Z$ is isomorphic to $Q \times^P   \u(\p)$.
Under this isomorphism, the projection from $Z$ to 
$Q \cdot \u(\p) = L \cdot \u(\p)$ 
is identified to the morphism $g'$, since it is $Q$-equivariant.
\end{proof}

We use notations  in Theorem \ref{Hi} (GP) in the following.
Assume furthermore that 
$\Delta_1 - \Gamma_1 \sim_R \Delta_1 - \Gamma'_1 $ is one of the 
fundamental relations in Theorem \ref{Hi}.
Let $\nu, \nu'$ be the  Springer maps associated to
 $P_{\Gamma_1 }$ and $P_{\Gamma'_1 }$ in the root system $\langle \Delta_1 \rangle  $.
Denote  $P = P_{\Gamma_1 \cup \Gamma_2}$,  
$P'= P_{\Gamma'_1 \cup \Gamma_2}$ and 
$Q =  P_{\Gamma_2}$.
Let $\pi, \pi'$ be the Springer maps associated to
$P, P'$ and $\0$ their Richardson orbit. 
\begin{Prop} \label{deg}
Under the above hypothesis, we have :

(i) $Q \cdot \u(P) = Q \cdot \u(P')$.

(ii) $\deg(\pi) \deg(\nu') = \deg(\pi') \deg(\nu)$.

(iii) The diagram $T^*(G/P) \xrightarrow{\phi} G \times^Q  (Q \cdot \u(P)) = 
 G \times^Q  (Q \cdot \u(P')) \xleftarrow{\phi'} T^*(G/P')$ is 
an analytically 
locally trivial family of the diagram given by $\nu$ and $\nu'$.
\end{Prop}
\begin{proof}
Let $Q = LU$ be a Levi decomposition of $Q$, then 
$U \subset P \cap P'$ since $P$ and $ P'$ are contained in  $ Q.$
 Notice that
$\p \cap \l(\q) $ and $\p' \cap \l(\q) $ are parabolic sub-algebras
in $\l(\q)$ corresponding to two polarizations of a same orbit,
so $L \cdot \u(\p \cap \l(\q))) = L \cdot \u(\p' \cap \l(\q))) .$ 
Now claim (i) follows from the fact $Q \cdot \u(P) = L \cdot \u(P)$
and $\u(P) =  \u(\p \cap \l(\q)) \oplus \u(\q)$.

Let $\psi:  G \times^Q  (Q \cdot \u(P)) \to \overline{\0}$
be the natural map. Then $\pi = \psi \circ \phi$ and  
$\pi' = \psi \circ \phi'$. 
Note that $\phi$ is the composition of the following maps:
$$
T^*(G/P) \simeq G \times^P \u(P) \simeq G \times^Q (Q \times^P \u(P))
\to G \times^Q (Q \cdot \u(P)),
$$
so the degree of $\phi$ is the same as that of the map
$Q \times^P \u(P) \to Q \cdot \u(P)$, which is equal to
the degree of the map 
$L \times^{(L \cap P)}  \u(\l(\q) \cap \p) \to  L \cdot
 (\u(\l(\q) \cap \p))$ by
Lemma \ref{lemdeg}.  The latter is in fact a trivial family
of the morphism $\nu$, which gives  $deg(\phi) = deg (\nu)$.
A similar argument shows that $deg(\phi') = deg (\nu')$,
 which gives (ii).
Note that the morphism $f$ in Lemma \ref{lemdeg}
is an affine bundle, so it is analytically locally trivial,
  which proves claim (iii).
\end{proof}
\begin{Rque}
One can prove directly assertion (ii) by using a formula of 
Borho-MacPherson (see Proposition \ref{Borho}).
\end{Rque}

This proposition is analogue to Proposition 4.4 
in \cite{Na}, but the proof is different here.
Now we turn to the proof of Theorem \ref{main}. 
We will argue case-by-case for the simple Lie algebra $\g$. 

If $\g$ is $A_n, G_2$ or $F_4$, then the relation $\sim$ 
is generated by (8) and (GP), by  Lemma \ref{equiv} and
the degree formula in Proposition \ref{deg}.

Assume $\g$ is $B_n$ and $P_\Gamma$ a polarization of $\0$,
where $$\Gamma =\{\alpha_{p_1}, \alpha_{p_1+p_2}, \cdots, 
\alpha_{p_1+\cdots+p_s} \}, p_i > 0, \forall i. $$ 
To simplify the notations, we encode $\Gamma$ by
the sequence of ordered numbers $[p_1, \cdots, p_s]$.
 If we want to perform the relation
$\sim_R$ in (1) of Theorem \ref{Hi}
 for  some $B_{3k-1}$, one should have 
$p_1+\cdots+p_{s-1} = n+1 - 3k$ and $p_s = 2k-1$ or $2k$.
We consider the case $p_s = 2k-1$, since the other one
can be done similarly. 
Then $P_{\Gamma^{(1)}} \sim_R P_\Gamma$, where 
$\Gamma^{(1)} =[p^{(1)}_1, \cdots, p^{(1)}_s] $ with $p^{(1)}_i = p_i$ for
$i \neq s$ and $p^{(1)}_s = 2k$.
Now any $P_{\Gamma^{(2)}}$ obtained from  $P_{\Gamma^{(1)}}$ by performing
(8) and (GP) in Theorem \ref{Hi} has the following form:
$\Gamma^{(2)} =[p^{(1)}_{\sigma(1)}, \cdots, p^{(1)}_{\sigma(s)}] $
for some element $\sigma$  in the symmetric group $\mathfrak{S}_s$. For simplicity,
we will denote by $\deg(\Gamma)$  the 
degree of the Springer map associated to $P_{\Gamma}$. Then we have 
$\deg(\Gamma^{(2)})=\deg(\Gamma^{(1)}) = 2 \deg(\Gamma)$. 

If we want to change the degree, we need to perform once again the 
operation in (1) of Theorem \ref{Hi}. 
There are only two possibilities: 
(i)  perform the operation in  (1) for  $B_{3k-1}$; (ii) perform 
the operation in (1) for $B_{3k-4}$.

For case (i), after the operation, the rightmost marked node goes back to its original
position and the degree remains the same.
For case (ii), after the operation, one obtains $\Gamma^{(3)}$, but the degree
goes higher: $\deg(\Gamma^{(3)})  = 2 \deg (\Gamma^{(2)}) = 4 \deg(\Gamma).$
By this way, we see that to obtain the 
same degree as $\deg(\Gamma)$ for $\Gamma^0 = [q_1, \cdots, q_s]$, one should have
$\sum_{j=1}^s q_j = \sum_{i=1}^s p_i$. In other words, the rightmost marked node should stay
at the same position. Now it follows that $\Gamma^0$ can be obtained from
$\Gamma$ by just performing operations in $(8)$ and $(GP)$ of
Theorem \ref{Hi}.

Similar arguments can be done to the case $\g = C_n$.
When $\g = D_n$,  we have two possible operations
(2) and (3) in Theorem \ref{Hi}
 which  do not preserve degrees. The key point is
 that if we have performed one of them, then we
can not perform the other one, so the situation is
similar to the $B_n$ case. 

If $\g$ is $E_6$, let $P_\Gamma$ be a parabolic subgroup.
The only possible
 operation not preserving the degree is (2) of Theorem \ref{Hi},
 since
(6) is settled by Lemma \ref{equiv}.
But then one should have $\alpha_1, \alpha_6 \in \Gamma$.
Now it is easy to see that for any $P_{\Gamma'}$ equivalent to
$P_{\Gamma}$, one can arrive $\Gamma'$ from $\Gamma$
by just performing  operations (8) and (GP) of Theorem \ref{Hi}. 

If $\g$ is $E_7$, we can perform either (2) or (6) of
 Theorem \ref{Hi}.
For (2), one should have $\alpha_1, \alpha_6 \in \Gamma$.
For (6), one has $\alpha_7 \in \Gamma$. The argument
is similar to the case of $E_6$.  We can do the similar
to the case of $E_8$, noticing that (7) is already done
by  Lemma \ref{equiv}. This finishes the proof of Theorem 
\ref{main}. \ \ \ \ \ Q.E.D. \vspace{0.3 cm}

Following Namikawa \cite{Na}, the diagrams given by the Springer maps
of  dual parabolic subgroups in (8), (9), (10) of Theorem \ref{Hi}
will be called {\em  stratified Mukai flops} of type $A, D,  E_{6, I}$ 
and $E_{6, II}$ respectively. The following is the list of the dual marked
Dynkin diagrams.

\[
\begin{minipage}{.35\textwidth}
    \begin{center}  
      \begin{diagram}[size=0.4cm]
\circ & \rLine & \cdots       &  \rLine        &  \bullet      &  \rLine  & \cdots & \rLine &\circ \\
 &      &   &  &   k &&&&
\end{diagram}
    \end{center}
  \end{minipage}
\begin{minipage}{.2\textwidth}
\begin{center}
   $A_{n-1,k} (2k \neq n)$
\end{center}
\end{minipage}
  \begin{minipage}{.35\textwidth}
    \begin{center}  
     \begin{diagram}[size=0.4cm]
     \circ & \rLine & \cdots       &      &  \bullet      &   & \cdots & \rLine &\circ \\
 &      &   &  &   n-k &&&&
\end{diagram}
    \end{center}
  \end{minipage}
\]

\[
\begin{minipage}{.35\textwidth}
    \begin{center}  
      \begin{diagram}[size=0.5cm]
\bullet      &       &        &         &       &   &      &  &\\
          &  \rdLine      &        &         &      & &         &  & \\
        &                 & \circ & \rLine & \circ &\cdots & \circ &\rLine  & \circ  \\
        &   \ruLine     &        &         &      & &          &&   \\
  \circ     &        &        &         &       &  &        &  & 
      \end{diagram}
    \end{center}
  \end{minipage}
\begin{minipage}{.15\textwidth}
\begin{center}
   $D_{2n+1}$
\end{center}
\end{minipage}
  \begin{minipage}{.35\textwidth}
    \begin{center}  
     \begin{diagram}[size=0.5cm]
      \circ      &       &        &         &  &     &     &    &  \\
          &  \rdLine      &        &         &  &     &     &    &   \\
        &                 & \circ & \rLine & \circ &\cdots & \circ & \rLine  & \circ  \\
        &   \ruLine     &        &         &       & &      &    &   \\
  \bullet     &        &        &         &       & &       &  &   
\end{diagram}
    \end{center}
  \end{minipage}
\]

\[
\begin{minipage}{.35\textwidth}
    \begin{center}  
      \begin{diagram}[size=0.5cm]
\bullet & \rLine & \circ & \rLine  & \circ & \rLine & \circ  & \rLine & \circ\\
        &        &       &         &  \vLine &      &        &        &   \\
        &        &       &         &  \circ &      &        &        & 
\end{diagram}
    \end{center}
  \end{minipage}
\begin{minipage}{.15\textwidth}
\begin{center}
   $E_{6, I}$
\end{center}
\end{minipage}
  \begin{minipage}{.35\textwidth}
    \begin{center}  
     \begin{diagram}[size=0.5cm]
\circ & \rLine & \circ & \rLine  & \circ & \rLine & \circ  & \rLine & \bullet\\
        &        &       &         &  \vLine &      &        &        &   \\
        &        &       &         &  \circ &      &        &        & 
\end{diagram}
    \end{center}
  \end{minipage}
\]

\[
\begin{minipage}{.35\textwidth}
    \begin{center}  
      \begin{diagram}[size=0.5cm]
\circ & \rLine & \bullet & \rLine  & \circ & \rLine & \circ  & \rLine & \circ\\
        &        &       &         &  \vLine &      &        &        &   \\
        &        &       &         &  \circ &      &        &        & 
\end{diagram}
    \end{center}
  \end{minipage}
\begin{minipage}{.15\textwidth}
\begin{center}
   $E_{6, II}$
\end{center}
\end{minipage}
  \begin{minipage}{.35\textwidth}
    \begin{center}  
     \begin{diagram}[size=0.5cm]
\circ & \rLine & \circ & \rLine  & \circ & \rLine & \bullet  & \rLine & \circ\\
        &        &       &         &  \vLine &      &        &        &   \\
        &        &       &         &  \circ &      &        &        & 
\end{diagram}
    \end{center}
  \end{minipage}
\]

In practice, the procedure to find all marked Dynkin diagrams
equivalent to a fixed one $\Gamma$ is the following: choose 
a node $\beta \in \Gamma$. Let $C$ be the maximal 
 connected subgraph containing $\beta$, with other nodes in
$\Delta - \Gamma$. Then $C$ is a single marked Dynkin diagram.
If $C$ is one of the above marked Dynkin diagram,
 we replace it  with the dual
one to obtain $\Gamma'$. Then we have $P_\Gamma \sim P_{\Gamma'}$ and
we can  continue the procedure with $\Gamma'$.

Let us deduce some corollaries.

\begin{Cor}\label{bir}
Assume that $\g$ is simple.
Let $\pi_i: T^*(G/P_i) \to \overline{\0}, i=1, 2$
 be two Springer maps
with the same degree, then we have  a  birational map 
$T^*(G/P_1) \dasharrow T^*(G/P_2)$ over $\overline{\0}$ which can be decomposed into
a finite sequence of analytically
 locally trivial families of stratified Mukai flops
of type $A, D,  E_{6, I}$ and $E_{6, II}$. 
\end{Cor}

This follows immediately from Theorem \ref{main}
and Proposition \ref{deg}. Note that for the special case
where $deg(\pi_i) = 1$, this implies that any two symplectic 
resolutions of a nilpotent orbit closure are connected by
stratified Muaki flops, which
has been previously proved in \cite{Na}. 
Our proof here is more elementary, in the sense that we do not
use Mori theory and the  Springer correspondence for exceptional Lie algebras.

Let $d$ be the degree of $\pi_i$ in the precedent lemma, then
$d$ divides the order of the fundamental group of  $\0$. 
Let $\0'$ be the $G$-covering
of degree $d$ of $\0$, which embeds into the unique open 
$G$-orbit in $T^*(G/P_i)$ (\cite{BK}).
The map  $\pi_i$ factorizes through
the symplectic resolution $T^*(G/P_i) \to \overline{\0'}$, 
where $\overline{\0'}$ is the image of the Stein factorization of $\pi_i$. 
 If $d$ is odd, then
every  symplectic resolution of $\overline{\0'}$
arises in this way (\cite{Fu2}), so in this case, the corollary implies that any two
symplectic resolutions of $\overline{\0'}$ are related by stratified Mukai flops.
\begin{Cor}\label{two}
Let $\0$ be a nilpotent orbit in a simple exceptional Lie algebra 
$\g$. Then the degrees of Springer maps associated to polarizations
of $\0$ can take at most two values.
\end{Cor}

This follows directly from the proof of Theorem \ref{main}. 
An immediately consequence is that $\overline{\0}$ has
at most two $G$-coverings which admit a Springer resolution.
Note that this corollary is not true if $\g$ is classical.

\begin{Cor}\label{exc}
Let $\0$ be a nilpotent orbit in a simple exceptional Lie algebra. Then $\overline{\0}$
admits a symplectic resolution if and only if $\0$ is Richardson and $\0$ is not
one of the following orbits:

$$ A_4+A_1, D_5(a_1) \text{ in }   E_7, \quad E_6(a_1)+A_1, E_7(a_3) 
\text{ in }   E_8.     $$
\end{Cor}
\begin{proof}

First consider the orbit $\0: = \0_{D_7(a_2)}$ in $E_8$.
A polarization of $\0$  is given by $P_{\alpha_1 \alpha_4}$
(\cite{Hi}). Let $\pi$ be the corresponding Springer map. 
Remark that we can perform the operation (3) in Theorem \ref{Hi} 
for $D_7$
to obtain another polarization $P_{\alpha_1 \alpha_5}$. 
If we denote by $\pi'$ the Springer map of $P_{\alpha_1 \alpha_5}$,
then by Proposition \ref{deg} and Lemma \ref{equiv},
we have $\deg(\pi) = 2 \deg(\pi')$. Notice that the fundamental
group of $\0$ is $S_2$ (\cite{CM}, p. 134), so the degree of 
any Springer map of $\0$ is either 1 or 2, 
which gives $deg(\pi') = 1$ and $\deg(\pi) =2$.
In particular, the closure of the orbit $\0_{D_7(a_2)}$ in $E_8$
admits a symplectic resolution.

Then consider the orbits in the  Corollary.
The following is a list of a polarization for each orbit
(\cite{Hi}):
\begin{center}
\begin{tabular}{|l|l|l|l|l|}
\hline 
algebra & 
 \multicolumn{2}{|l|}{$E_7$}  & \multicolumn{2}{l|}{$E_8$}  \\
\hline
orbit & $A_4+A_1$  & $D_5(a_1)$  & $E_6(a_1)+A_1$ &
   $E_7(a_3)$   \\
\hline
polarization & $ P_{\alpha_2 \alpha_3}$ & $P_{\alpha_1 \alpha_2 \alpha_3}$
& $P_{\alpha_1 \alpha_2 \alpha_4}$ & 
 $P_{\alpha_1 \alpha_2 \alpha_3 \alpha_4}$\\
\hline
\end{tabular} 
\end{center}
Now we do a case-by-case check to show that
the relation $\sim_R$ coincides with $\sim $, i. e.
we can not perform the operations  
(2), (3) (for $D_7$), (6) and  (7)
in Theorem \ref{Hi} to any polarization of the orbit.
This implies that for any two polarizations $Q_1, Q_2$  of 
one of the four orbits, the degrees of the associated Springer
maps are the same.
In the appendix, we calculate the degrees of the Springer maps
associated to the above polarizations
by using a formula in \cite{BM}, which turn out to be 2. 
Thus the four
orbit closures do not admit a symplectic resolution.

The proof is completed by applying Proposition \ref{exc_f} and
Theorem \ref{exceptional}.
\end{proof}

\begin{Exam} \label{D7} \upshape
Let  $\0 = \0_{D_7(a_2)}$ in $E_8$ and $\0_s$ its universal
$G$-covering. Then $\overline{\0}$
admits exactly two symplectic resolutions, given by
 $P_{\alpha_1 \alpha_5}$ and $P_{\alpha_2 \alpha_5}$
(with marked Dynkin diagrams listed in the following),
which is  a locally trivial family of 
ordinary Mukai flops of $T^*\pit^4$. Notice that 
neither of the two parabolic sub-groups contains 
the Jacobson-Morozov parabolic sub-group $P$ of $\0$, so
we cannot obtain a symplectic resolution of $\overline{\0}$
by extremal contractions of $\mu: G \times^P \n \to \overline{\0}.$

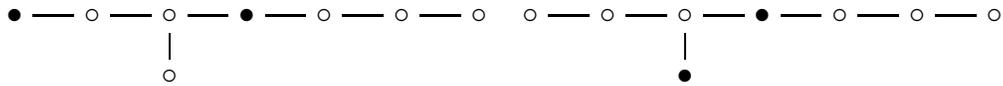
\[
\begin{minipage}{.35\textwidth}
    \begin{center}  
      \begin{diagram}[size=0.4cm]
\bullet & \rLine & \circ & \rLine  & \circ & \rLine & \bullet  & \rLine & \circ & \rLine & \circ & \rLine & \circ\\
        &        &       &         &  \vLine &      &        &        &  &  & && \\
        &        &       &         &  \circ &      &        &        & 
\end{diagram}
    \end{center}
  \end{minipage}
\begin{minipage}{.15\textwidth}
\begin{center}
\end{center}
\end{minipage}
  \begin{minipage}{.35\textwidth}
    \begin{center}  
     \begin{diagram}[size=0.4cm]
\circ & \rLine & \circ & \rLine  & \circ & \rLine & \bullet  & \rLine & \circ & \rLine & \circ & \rLine & \circ\\
        &        &       &         &  \vLine &      &        &        &  &  & && \\
        &        &       &         &  \bullet &      &        &        & 
\end{diagram}
    \end{center}
  \end{minipage}
\]

The closure $\overline{\0}_s$ admits exactly four
different Springer resolutions, given by 
$P_{\alpha_1 \alpha_4}$, $P_{\alpha_3 \alpha_4}$,
 $P_{\alpha_3 \alpha_7}$ and $P_{\alpha_5 \alpha_7}$
(their marked Dynkin diagrams are listed in the following),
where the diagram of two consequentive 
symplectic resolutions is 
a locally trivial family of stratified Mukai flops of $T^*(\pit^2),
T^*(Gr(2, 7))$ and of type $E_{6, II}$ respectively.
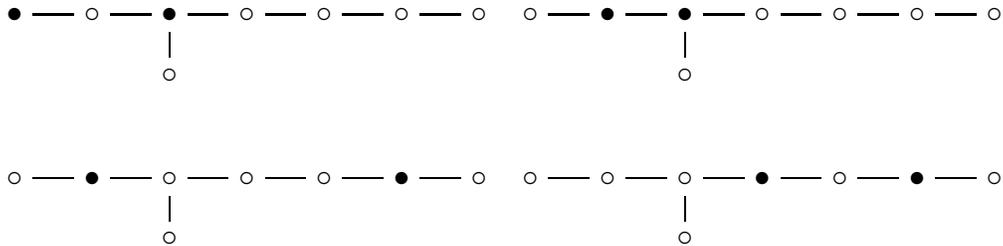
\[
\begin{minipage}{.35\textwidth}
    \begin{center}  
      \begin{diagram}[size=0.4cm]
\bullet & \rLine & \circ & \rLine  & \bullet & \rLine & \circ  & \rLine & \circ & \rLine & \circ & \rLine & \circ\\
        &        &       &         &  \vLine &      &        &        &  &  & && \\
        &        &       &         &  \circ &      &        &        & 
\end{diagram}
    \end{center}
  \end{minipage}
\begin{minipage}{.15\textwidth}
\begin{center}
\end{center}
\end{minipage}
  \begin{minipage}{.35\textwidth}
    \begin{center}  
     \begin{diagram}[size=0.4cm]
\circ & \rLine & \bullet & \rLine  & \bullet & \rLine & \circ  & \rLine & \circ & \rLine & \circ & \rLine & \circ\\
        &        &       &         &  \vLine &      &        &        &  &  & && \\
        &        &       &         &  \circ &      &        &        & 
\end{diagram}
    \end{center}
  \end{minipage}
\]

\[
\begin{minipage}{.35\textwidth}
    \begin{center}  
      \begin{diagram}[size=0.4cm]
\circ & \rLine & \bullet & \rLine  & \circ & \rLine & \circ  & \rLine & \circ & \rLine & \bullet & \rLine & \circ\\
        &        &       &         &  \vLine &      &        &        &  &  & && \\
        &        &       &         &  \circ &      &        &        & 
\end{diagram}
    \end{center}
  \end{minipage}
\begin{minipage}{.15\textwidth}
\begin{center}
\end{center}
\end{minipage}
  \begin{minipage}{.35\textwidth}
    \begin{center}  
     \begin{diagram}[size=0.4cm]
\circ & \rLine & \circ & \rLine  & \circ & \rLine & \bullet  & \rLine & \circ & \rLine & \bullet & \rLine & \circ\\
        &        &       &         &  \vLine &      &        &        &  &  & && \\
        &        &       &         &  \circ &      &        &        & 
\end{diagram}
    \end{center}
  \end{minipage}
\]

\end{Exam}

\section{A description of stratified Mukai flops}

By Corollary \ref{bir}, to understand the birational geometry
of Springer maps of the same degree, one is led to understand
 stratified Mukai flops of type $A, D, E_{6, I}$ and $ E_{6, II}$. We will give
a uniform resolution and explicit description of these flops. Note that some other 
descriptions of these flops were presented in \cite{Ch}. 

Let $Q = P_{\alpha_i}, Q' = P_{\alpha_j}$
 be a pair of the dual standard parabolic subgroups in (8), (9), (10)
of Theorem \ref{Hi} with Lie algebras $\q, \q'$
 and $\0$ their
Richardson orbit.
We denote by $\pi$ and $ \pi'$ the Springer maps 
associated to $Q$ and $Q'$.
 Write $P$ the Jacobson-Morozov parabolic 
subgroup of $\0$ with Lie algebra $\p$
 and $\n$ the natural nilpotent ideal in
$\u(P)$. 

\begin{Thm}\label{Mukai}
(i) The weighted Dynkin diagram of $\0$ has weight 
$1$ on nodes $\alpha_i, \alpha_j$ and $0$ on other nodes.

(ii) $\q$ and $\q'$ are the only (non-trivial) 
standard parabolic sub-algebras containing
$\p$  and $\n = \u(\q) \cap \u(\q').$

(iii) The variety $G \times^P \n$ is isomorphic to 
the graph closure of the flop
 $\phi: T^*(G/Q) \dasharrow T^*(G/Q') $.

(iv) We have the following commutative diagram:
\begin{diagram}[size=0.6cm]
& &  &   G  \times^P \n &  &  &  \\
&&  \ldTo^{\nu} &      &  \rdTo^{\nu'}  &&  \\
 G \times^Q (Q \cdot \n) \simeq &T^*(G/Q)  & & \dTo_{\mu}  & & 
T^*(G/Q')& \simeq G \times^{Q'}(Q' \cdot \n) \\
& & \rdTo_{\pi} &      &  \ldTo_{\pi'}  &&  \\
& &  &    \overline{\0} &    & & \\
\end{diagram}
The natural morphisms $\nu, \nu'$ are contractions of 
extremal rays of $\mu$ and they are also the two graph
 projections under the isomorphism in (iii).
\end{Thm}
\begin{proof}
The following is a list of Richardson orbits appeared
in the stratified Mukai flops.
\begin{center}
\begin{tabular}{|l|l|l|l|l|}
\hline 
type & $ A_{n-1, k} (k < n/2) $ 
 &  $D_{2n+1}$  &  $E_{6, I}$  &  $E_{6, II}$ \\
\hline
orbit &  $\0_{[2^k, 1^{n-2k}]}$ &  $\0_{[2^{2n}, 1^2]}$ 
& $\0_{2A_1}$ & $\0_{A_2+2A_1}$ \\
 \hline
\end{tabular} 
\end{center}
Now claim (i) follows immediately(see \cite{CM}). The first part of 
claim (ii) follows from (i). 
Claim (iv) and the inclusion 
$\n \subseteq \u(Q) \cap \u(Q')$
have already been verified in Theorem \ref{classical}
and Theorem \ref{exceptional} while the  inclusion 
$ \u(Q) \cap \u(Q') \subseteq \n$ is obvious.
  
To show (iii), note that the two projections $G/P \to G/Q, G/P \to G/Q$ embed 
$G/P$ as a sub-variety in $G/Q \times G/Q'$ (the incident variety).
The composition of morphisms  
$$G \times^P \n \to G/P \times \overline{\0}
\xrightarrow{\eta} (G/Q \times \overline{\0}) \times (G/Q' \times \overline{\0}) $$ 
embeds $G \times^P \n$ into a closed sub-variety of 
$T^*(G/Q) \times T^*(G/Q')$ (here one uses (ii)), where 
$\eta$ is given by $\eta([gP], x) = (([gQ], x), ([gQ'], x))$.
The image is in fact  the graph closure of the flop $\phi$
in  
$T^*(G/Q) \times_{\overline{\0}} T^*(G/Q')$, since it
is irreducible closed and 
contains the diagonal embedding of $\0$ into the fiber product.

\end{proof}

\begin{Rque} \upshape
(i) It is possible to show that stratified Mukai flops
 are the only flops which
appear in symplectic resolutions of nilpotent orbit closures 
with properties (i), (ii), (iv) in Theorem \ref{Mukai}.

(ii) A similar diagram holds for some other flops (see Theorem
\ref{classical} and  \ref{exceptional}).

(iii)  The variety $G \times^P \n$ is in fact isomorphic to the conormal 
bundle of $G/P$ in $G/Q \times G/Q'$.
\end{Rque}

\begin{Exam} \upshape
Consider the Mukai flop of type  $A_n$ 
with $i=1, j=n$. Then $G \times^P \n$
is isomorphic to  the blowup of $T^*(\pit^n)$ along the zero section and
$G/P$ is the incidence variety in $\pit^n \times (\pit^n)^*$. 
The two extremal contractions are blow-downs of $G/P$ along two different directions.

Let $q: T^*(G/Q) \to G/Q $ and $q': T^*(G/Q') \to G/Q'$ be the two projections and
$\Phi= \nu'_* \nu^*: K(T^*(G/Q)) \to K(T^*(G/Q'))$
 the natural morphism between two $K$-groups. It has been observed in
\cite{Na1} that $\Phi$ is not an isomorphism even when $n=2$.
In fact,  for $-n+1 \leq k \leq 0$, we have $\Phi(q^* \0(k)) = (q')^* \0(-k)$ and 
$\Phi(q^* \0(1)) = (q')^* \0(-1) \otimes \mathcal{I}_{G/Q'}$, where $\mathcal{I}_{G/Q'}$ is the ideal sheaf of 
$G/Q'$ in $T^*(G/Q')$. Now using the Koszul resolution and the Euler exact sequence, we can obtain that
$\Phi(q^* \0(1)) = - n  (q')^* \0(-1) + \sum_{i=0}^{n-1} k_i (q')^*(\0(i))$ for some integers $k_i$, which shows that
$\Phi$ is not an isomorphism as soon as $n \geq 2$. However, the morphism 
$$ \Phi \otimes_\zit \qit:  K(T^*(G/Q)) \otimes_{\zit} \qit \to K(T^*(G/Q')) \otimes_{\zit} \qit $$
is an isomorphism. Similarly for the $G$-equivariant $K$-groups, we have also an isomorphism: 
$$\Phi_G \otimes_\zit \qit: K^G(T^*(G/Q)) \otimes_{\zit} \qit \to K^G(T^*(G/Q')) \otimes_{\zit} \qit. $$ 
It seems plausible that similar results hold for other stratified Mukai flops.
More generally, for two birational $K$-equivalent varieties, it seems that the graph closure gives an isomorphism
between the $K$-groups with rational coefficients
(see \cite{Wan} for further discussions and related references). 
\end{Exam}

\section{Appendix: The degrees of Springer maps}
Let $W$ be the Weyl group of $G$.
The Springer correspondence (\cite{Sp2}) assigns to
any irreducible $W$-module a unique  pair  $(\0, \phi)$ 
consisting of a nilpotent orbit 
$\0$ in $\g$ and an irreducible representation $\phi$ of the component group $A(\0):= G^x/(G^x)^\circ$ of $\0$,
where  $x$ is any point in $\0$ and $(G^x)^\circ$ is the identity component of $G^x$.
 The corresponding irreducible $W$-module will be denoted by  $\rho_{(\0, \phi)}.$ This correspondence is not surjective
onto the set of all pairs $(\0, \phi)$. A pair will be called {\em relevant} if it corresponds 
to an irreducible $W$-module, then
the Springer correspondence establishes a bijection between
irreducible $W$-modules and relevant pairs in $\g$.
For exceptional cases, the Springer correspondence
has been completely worked out in \cite{Sp1} for $G_2$, in \cite{Sho} for $F_4$ and in \cite{AL} for $E_n (n=6,7,8)$. 

Consider a parabolic sub-group $Q$ in $G$.
Let $L$ be a Levi sub-group of $Q$ and $T$ a maximal torus in $L$.  The Weyl group of 
$L$ is $W(L) := N_L(T)/T$, where $N_L(T)$ is the normalizer of $T$ in $L$. It is a sub-group
of the Weyl group $W$ of $G$. Let $\epsilon_{W(L)}$ be the sign representation of $W(L)$
and $\Ind_{W(L)}^W(\epsilon_{W(L)})$ the induced representation of $\epsilon_{W(L)}$ to $W$.

\begin{Prop}[\cite{BM}, Corollary 3.9]\label{Borho}
Let $\pi_Q: T^*(G/Q) \to \overline{\0_Q}$ be 
the Springer map associated to the parabolic sub-group $Q$.
Then  $$\deg(\pi_Q) = \sum_{\phi} \mtp(\rho_{(\0_Q, \phi)}, 
\Ind_{W(L)}^W(\epsilon_{W(L)})) \deg \phi,$$
where the sum is over all irreducible representations $\phi$ of $A(\0_Q)$ such that $(\0_Q, \phi)$   
is a relevant pair, 
$\mtp(\rho_{(\0_Q, \phi)},\Ind_{W(L)}^W(\epsilon_{W(L)})) $
 is the multiplicity
 of $ \rho_{(\0_Q, \phi)} $ 
in $ \Ind_{W(L)}^W(\epsilon_{W(L)})$ 
and $\deg \phi$ is the dimension of 
 the irreducible representation $\phi$.
\end{Prop}

The multiplicity $ \mtp(\rho_{(\0_Q, \phi)}, \Ind_{W_0}^W(\rho))$ has been worked out in \cite{Alv}, for
any  irreducible representation $\rho$ of any maximal parabolic sub-group  $W_0$ of $W$, where
$ \Ind_{W_0}^W(\rho)$ is the induced representation of $\rho$ to $W$. Note that 
$ \Ind_{W(L)}^W(\epsilon_{W(L)}) = \Ind_{W_0}^W (\Ind_{W(L)}^{W_0}(\epsilon_{W(L)}))$ for any maximal parabolic sub-group
$W_0$ of $W$ containing $W(L)$ and
 $\Ind_{W(L)}^{W_0}(\epsilon_{W(L)}) $ can be determined by the Littlewood-Richardson rules 
when $W_0$ is classical and by \cite{Alv} when $W_0$ is exceptional.
 Using this,  we can calculate the degrees of the Springer maps
associated to the parabolic subgroups appeared in Theorem \ref{Hi}, 
and the result is as follows:

\begin{center}
\begin{tabular}{|l|l|l|l|l|l|}
\hline 
 \text{Lie algebra} & $G_2$ &  $F_4$   &  $F_4$  & $E_6$ & $E_8$  \\
\hline
\text{parabolic subgroup}  &  $P_{\alpha_1}$ &  $P_{\alpha_3}$ & $P_{\alpha_1 \alpha_4}$  &  $P_{\alpha_2 \alpha_5}$ &   $P_{\alpha_2 \alpha_3}$   \\
\hline
\text{nilpotent orbit}  & $G_2(a_1) $ & $F_4(a_3)$   & $ F_4(a_3)$ & $D_4(a_1)$ & $E_8(a_7)$ \\
 \hline
\text{component group} &  $S_3$ & $S_4$ & $S_4$ & $S_3$ & $S_5$ \\
\hline
\text{degree}  &  $2$  &  $ 4 $  &  $6$  &  $3$   &  $10$ \\
\hline
\end{tabular} 
\end{center}

In a similar way, we obtain the degrees of
the Springer maps for the orbits in the list
of Corollary \ref{exc}:

\begin{center}
\begin{tabular}{|l|l|l|l|l|}
\hline 
 \text{Lie algebra} & $E_7$ &  $E_7$   &  $E_8$  & $E_8$  \\
\hline
\text{parabolic subgroup}  &  $P_{\alpha_2 \alpha_3}$  &  $P_{\alpha_1 \alpha_2 \alpha_3} $ & $P_{\alpha_1 \alpha_2 \alpha_4}$  &  $P_{\alpha_1 \alpha_2 \alpha_3 \alpha_4}$    \\
\hline
\text{nilpotent orbit}  & $A_4+A_1 $ & $D_5(a_1)$   & $ E_6(a_1)+A_1$ & $E_7(a_3)$  \\     
\hline
\text{degree}  &  $2$  &  $ 2 $  &  $2$  &  $2$    \\
\hline
\end{tabular} 
\end{center}
\begin{Rque} \upshape
The correspondence between notations of irreducible characters of $E_n (n=6,7,8)$ in \cite{Alv} and
those in \cite{AL} is given in \cite{BL}. See also \cite{GP} (Appendix $C$).
\end{Rque}

\begin{Rque}\label{det} \upshape
When  $\g$ is classical, the degree of a Springer map is
given by  Hesselink's formula (\cite{He}), in terms of the flag type of the parabolic
sub-group. Our method here allows one to find the degree of the Springer map from the marked Dynkin diagram of
the parabolic sub-group in exceptional Lie algebras.
 (This also works  for
Lie algebras of classical type).
In fact, when $\g$ is exceptional, for any Richardson
orbit, we have either given  a degree one polarization 
or proved the degree of any polarization is 2 (Theorem \ref{exceptional} and 
Corollary \ref{exc}). Now we can use Proposition \ref{deg},
 Lemma \ref{equiv} and the above results 
 to determine the degree associated to any other polarization.
\end{Rque}
\begin{Exam} \upshape
We will calculate the degree $d$ of the Springer map
associated to the following parabolic sub-group $Q$ in $E_7$: 
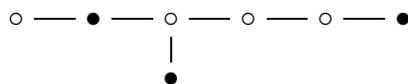
\[
     \begin{diagram}[size=0.4cm]
\circ & \rLine  & \bullet & \rLine & \circ  & \rLine & \circ & \rLine & \circ & \rLine &\bullet\\
        &        &       &         &  \vLine &      &        &        &  &  & && \\
        &        &       &         &  \bullet &      &        &        & 
\end{diagram}
\]

We can perform the operation (6) in Theorem \ref{Hi} to obtain the 
following parabolic $Q'$ in $E_7$:
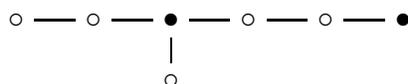
\[  
     \begin{diagram}[size=0.4cm]
 \circ & \rLine  & \circ & \rLine & \bullet & \rLine & \circ & \rLine & \circ & \rLine & \bullet\\
        &        &       &         &  \vLine &      &        &        &  &  & && \\
        &        &       &         &  \circ &      &        &        & 
\end{diagram}
\]

Note that $Q'$ is the Jacobson-Morozov parabolic sub-group of 
the even orbit $\0:=E_7(a_5)$, so its associated Springer map is 
birational. By Proposition \ref{deg}, we get that the degree of the Springer
map associated to $Q$ is  $d =3$. 
The component group $A(\0)$ is isomorphic to $S_3$, so the 
$2$-fold, $6$-fold $G$-coverings of $\overline{\0}$
have no Springer resolution  by Corollary \ref{two}.

\end{Exam}

\quad \\
C.N.R.S.,
Laboratoire J. Leray (Math\'ematiques)\\
 Facult\'e  des sciences, Univ. de Nantes \\
2, Rue de la Houssini\`ere,  BP 92208 \\
F-44322 Nantes Cedex 03 - France\\
\quad \\
fu@math.univ-nantes.fr
\end{document}